\documentclass[11pt,a4paper,reqno]{amsart}

\usepackage[utf8]{inputenc}
\usepackage{amsmath,amssymb,amsthm,mathtools}
\usepackage{enumitem}
\usepackage[a4paper,margin=3cm]{geometry}
\usepackage{microtype}
\usepackage[colorlinks=true,linkcolor=blue,citecolor=blue,urlcolor=blue]{hyperref}
\usepackage[nameinlink,capitalise]{cleveref}

\theoremstyle{plain}
\newtheorem{theorem}{Theorem}[section]
\newtheorem{lemma}[theorem]{Lemma}
\newtheorem{proposition}[theorem]{Proposition}
\newtheorem{corollary}[theorem]{Corollary}
\newtheorem*{conjecture}{Conjecture}
\theoremstyle{definition}

\newtheorem{remark}[theorem]{Remark}

\numberwithin{equation}{section}
\setlist[enumerate,1]{label=\textup{(\arabic*)},leftmargin=*}

\newcommand{\R}{\mathbb{R}}
\newcommand{\N}{\mathbb{N}}

\newcommand{\II}{\mathbb{I}}
\newcommand{\PP}{\mathbb{P}}
\newcommand{\EE}{\mathbb{E}}
\newcommand{\ind}{\mathbf{1}}
\newcommand{\dd}{\mathrm{d}}
\newcommand{\cB}{\mathcal{B}}
\newcommand{\cP}{\mathcal{P}}

\newcommand{\cW}{\mathcal{W}}
\newcommand{\Cd}{\mathcal{C}^{d}}
\newcommand{\Car}{\mathcal{C}^{d}_{ar}}
\newcommand{\FK}{F_K^{d}}
\newcommand{\bx}{\mathbf{x}}
\newcommand{\bU}{\mathbf{U}}
\newcommand{\Pk}[1]{\Psi^{[#1]}}

\begin{document}

\title[Kendall distribution functions and Archimedean copulas]{Archimedean
copulas are uniquely determined by their Kendall distribution function in
dimension $d \geq 4$}

\author{Nicolas Pascal Dietrich}
\address{Jakob-Haringer-Straße 6, Techno 6, 5020 Salzburg, Austria}
\email{nicolaspascal.dietrich@plus.ac.at}

\begin{abstract}
Motivated by the well-known conjecture saying that every $d$-dimensional
Archimedean copula is uniquely determined by its Kendall distribution
function -- a property which so far has only been established for $d=2$ and
$d=3$ -- this contribution proves the conjecture for arbitrary dimension
$d\geq 2$. Working with Williamson measures instead of generators we show
that, given two Williamson measures inducing the same Kendall distribution
function, the scaling factors matching their Williamson-type transforms of
order $k\in\{1,\ldots,d-1\}$ form a cyclic pursuit system. Monotonicity
properties of this system directly yield the result in the non-strict case.
In the strict case a priori bounds in combination with a lower bound for the
pursuit rates are then shown to imply that the system is in consensus, i.e.,
that the two Williamson measures coincide up to scaling.
\end{abstract}


\maketitle

\setcounter{tocdepth}{1}

\section{Introduction}\label{sec:intro}

Archimedean copulas form one of the most frequently used classes of
dependence models, which, not least, goes back to the fact that they are easy
to handle analytically: every $d$-dimensional Archimedean copula $C$ can be
expressed in terms of a single univariate function
$\psi\colon[0,\infty)\to[0,1]=:\II$, its so-called generator, via
\[
  C(x_1,\ldots,x_d)=\psi\big(\varphi(x_1)+\cdots+\varphi(x_d)\big),
\]
whereby $\varphi$ denotes the pseudo-inverse of $\psi$. Hence, many questions
concerning Archimedean copulas boil down to questions about univariate
objects (see, e.g., \cite{KDT,MN,Nelsen}). One such univariate object is the
Kendall distribution function $\FK$ of a copula $C$, i.e., the distribution
function of the random variable $C(\bU)$, where $\bU$ is a random vector with
distribution function $C$. Kendall distribution functions are a key tool for
inference in Archimedean copula models (see \cite{GNZ,GR}), which naturally
triggers the question whether a $d$-dimensional Archimedean copula is
uniquely determined by its Kendall distribution function. Genest and Rivest
\cite{GR} showed that in the bivariate setting this is indeed the case.
Genest, Ne\v{s}lehov\'a and Ziegel \cite{GNZ} proved the same property in the
trivariate setting and provided evidence that it holds in arbitrary
dimension, which led to the following conjecture (also see
\cite[Remark~5.6]{KDT}), whereby $\Car$ denotes the family of all
$d$-dimensional Archimedean copulas:

\begin{conjecture}
Let $d\geq 4$ and $C_1,C_2\in\Car$. If $C_1$ and $C_2$ have the same Kendall
distribution function, then $C_1=C_2$.
\end{conjecture}

To the best of our knowledge the conjecture has been open for every
$d\geq 4$. In \cite[Remark~5.6]{KDT} it was pointed out that the conjecture
would follow from injectivity of the mapping assigning each (normalized)
Williamson measure $\gamma$ the Kendall distribution function of the
Archimedean copula induced by $\gamma$. Building upon this very idea, in the
current contribution we prove the conjecture in full generality:

\begin{theorem}\label{thm:main}
Let $d\geq 2$ and $C_1,C_2\in\Car$. If $C_1$ and $C_2$ have the same Kendall
distribution function, then $C_1=C_2$.
\end{theorem}

In terms of Williamson measures \cref{thm:main} states that two probability
measures $\gamma_1,\gamma_2$ on $(0,\infty)$ induce the same Kendall
distribution function if, and only if, they coincide up to scaling (see
\cref{thm:mainW} and \cref{cor:injective}).

The proof can be sketched as follows: Let $\gamma_1,\gamma_2$ be Williamson
measures of $C_1$ and $C_2$ with distribution functions $G_1,G_2$, let $K$
denote the common Kendall distribution function, set $n:=d-1$ and consider
the transforms $\Pk{k}_i(s):=\int(1-t/s)_+^{k}\,\dd\gamma_i(t)$ for
$k\in\{1,\ldots,n\}$; notice that, up to the substitution $s\mapsto 1/s$, the
function $\Pk{n}_i$ is the generator of $C_i$. Our starting point is the
identity $G_i=K\circ\Pk{n}_i$ (\cref{cor:key}). For every $k$ and every $s$ we
let $\varrho_k(s)$ denote the scaling factor fulfilling
$\Pk{k}_1(s)=\Pk{k}_2(\varrho_k(s)\,s)$, so that $\gamma_1$ and $\gamma_2$
coincide up to scaling if, and only if, all functions $\varrho_k$ are equal to
the same constant. Combining the afore-mentioned identity with the
differentiation rule
$s\,\frac{\dd}{\dd s}\Pk{k}_i=k\,(\Pk{k-1}_i-\Pk{k}_i)$ we show that the
(right-hand) derivative of $\varrho_k$ is non-negative at every point at which
$\varrho_{k-1}\geq\varrho_k$ and non-positive at every point at which
$\varrho_{k-1}\leq\varrho_k$, whereby $\varrho_0:=\varrho_n$. Loosely
speaking, each factor pursues its predecessor and the first factor pursues
the last one (\cref{prop:pursuit}). As a consequence, $\max_k\varrho_k$ is
non-increasing and $\min_k\varrho_k$ is non-decreasing (\cref{cor:monotone}).
If one of the two copulas is non-strict, then all factors start from the same
value at the lower end of the support, and monotonicity immediately yields
the result (\cref{prop:nonstrict}). In the strict case the pursuit system
lives on the whole of $(0,\infty)$ and an additional argument is needed: we
first derive a priori bounds showing that
$\max_k\varrho_k\leq n^2\min_k\varrho_k$ (\cref{lem:apriori}), then prove that
for $k\geq 2$ the factor $\varrho_k$ approaches $\varrho_{k-1}$ from above at
a rate which is bounded from below (\cref{lem:rate}), and finally use a
cascade argument (\cref{lem:cascade}) to show that $\varrho_1(s)$ converges to
$\sup_{t,k}\varrho_k(t)$ for $s\downarrow 0$ (\cref{prop:top}). Interchanging
the roles of the two measures it follows that $\varrho_1(s)$ converges to
$\inf_{t,k}\varrho_k(t)$ too, so all factors are equal to the same constant
(\cref{thm:strict}).

The rest of this contribution is organized as follows: \cref{sec:prelim}
gathers notation and preliminaries on Williamson measures that are used
throughout the text. In \cref{sec:pursuit} we introduce the afore-mentioned
pursuit system and derive its monotonicity properties. \cref{sec:proof} then
tackles the non-strict and the strict case separately and contains the proof
of \cref{thm:main}.

\section{Notation and preliminaries}\label{sec:prelim}

\subsection{Archimedean copulas and Kendall distribution functions}

Throughout the whole contribution $d\geq 2$ denotes a fixed integer and, in
order to keep notation as simple as possible, we set
\[
  n:=d-1 .
\]
As usual, $\Cd$ denotes the family of all $d$-dimensional copulas and vectors
are written in bold symbols. For every $C\in\Cd$ the corresponding
$d$-stochastic measure will be denoted by $\mu_C$. Furthermore we write
$\II:=[0,1]$, let $\lambda$ denote the Lebesgue measure and, for every metric
space $S$, $\cB(S)$ the Borel $\sigma$-field on $S$. For a function $f$ of
one real variable $f(s^-)$ denotes the left-hand limit and $f'_+(s)$ the
right-hand derivative of $f$ at $s$ (whenever they exist), and
$x_+:=\max\{x,0\}$. For more background on copulas we refer to
\cite{DS,Nelsen}.

Following \cite{KDT,MN}, a function $\psi\colon[0,\infty)\to\II$ will be
called (Archimedean) generator if it is continuous and non-increasing,
fulfills $\psi(0)=1$ as well as $\lim_{z\to\infty}\psi(z)=0=:\psi(\infty)$,
and is strictly decreasing on the interval
$[0,\inf\{z\in[0,\infty]\colon \psi(z)=0\})$. For every generator $\psi$ its
pseudo-inverse $\varphi\colon\II\to[0,\infty]$ is defined by
$\varphi(y):=\inf\{z\in[0,\infty]\colon\psi(z)=y\}$. We call $C\in\Cd$
Archimedean, and write $C\in\Car$, if some generator $\psi$ exists such that
$C=C_\psi$, whereby
\[
  C_\psi(\bx):=\psi\big(\varphi(x_1)+\cdots+\varphi(x_d)\big),
  \qquad \bx\in\II^d .
\]
In case of $\varphi(0)=\infty$ we refer to $C$ as strict and as non-strict
otherwise. The Kendall distribution function $\FK\colon\II\to\II$ of
$C\in\Car$ is defined by
\[
  \FK(t):=\mu_C\big(\{\bx\in\II^d\colon C(\bx)\leq t\}\big),
\]
i.e., $\FK$ is the distribution function of the random variable $C(\bU)$ for
$\bU\sim C$.

\subsection{Williamson measures}

Following \cite{Williamson} (also see \cite{KDT,MN}), generators of
$d$-dimensional Archimedean copulas can be described and handled in terms of
probability measures. In the sequel we will let $\cP$ denote the family of
all probability measures on $\cB((0,\infty))$. Notice that $\cP$ can be
identified with the family of all probability measures $\gamma$ on
$\cB([0,\infty))$ fulfilling $\gamma(\{0\})=0$, which is the formulation used
in \cite{KDT} -- in what follows we will use this identification without
further mention. The subsequent characterization goes back to
\cite{Williamson} and \cite{MN} (see \cite[Proposition~3.1 and
Corollary~3.1]{MN}, where the result is formulated in terms of the
distribution of $1/T$ for $T\sim\gamma$, as well as
\cite[Theorem~5.1]{KDT}):

\begin{theorem}\label{thm:williamson}
A function $\psi\colon[0,\infty)\to\II$ is the generator of some
$C\in\Car$ if, and only if, there exists some $\gamma\in\cP$ $($i.e., some
probability measure $\gamma$ on $[0,\infty)$ with $\gamma(\{0\})=0)$ such that
\begin{equation}\label{eq:williamson}
  \psi(z)=\int_{(0,\infty)}(1-tz)_+^{\,n}\,\dd\gamma(t)=:(\cW_d\gamma)(z)
\end{equation}
holds for every $z\geq 0$. In this case the measure $\gamma$ is unique.
\end{theorem}

Notice that the requirement $\gamma(\{0\})=0$, which is built into the
definition of $\cP$, cannot be dropped: if $\gamma$ is a probability measure
on $[0,\infty)$ with $\gamma(\{0\})>0$, then
$\int_{[0,\infty)}(1-tz)_+^{\,n}\,\dd\gamma(t)$ converges to
$\gamma(\{0\})>0$ for $z\to\infty$ and hence does not define a generator. In
the sequel this property will be used in the form
$\lim_{s\downarrow 0}\gamma((0,s])=0$ as well as in the form that every
random variable $T$ with distribution $\gamma\in\cP$ fulfills $T>0$ almost
surely.

For every $\gamma\in\cP$ we will write $\psi_\gamma:=\cW_d\gamma$ as well as
$C_\gamma:=C_{\psi_\gamma}\in\Car$, and refer to $\gamma$ as \emph{Williamson
measure} of $C_\gamma$. Contrary to \cite{KDT} we do not normalize
generators, so Williamson measures in our sense are only unique up to
scaling (see \cref{lem:scaling}). Furthermore, for $\gamma\in\cP$, $s>0$ and
$k\in\N$ we set
\begin{align*}
  G_\gamma(s)&:=\gamma((0,s]), &
  s_\gamma&:=\inf\{s>0\colon G_\gamma(s)>0\}\in[0,\infty),\\
  \Pk{k}_\gamma(s)&:=\int_{(0,s)}\Big(1-\frac{t}{s}\Big)^{k}\,\dd\gamma(t), &
  \Psi_\gamma&:=\Pk{n}_\gamma ,
\end{align*}
so that $s_\gamma$ is the left endpoint of the support of $\gamma$ and
$\Psi_\gamma(s)=\psi_\gamma(1/s)$ holds for every $s>0$. Notice that
$s_\gamma=0$ only means that $\gamma$ has mass arbitrarily close to $0$ --
mass at the point $0$ itself is excluded for every $\gamma\in\cP$. If $T$ is
a random variable with distribution $\gamma$ (for short: $T\sim\gamma$), then
obviously $\Pk{k}_\gamma(s)=\EE\big[(1-T/s)_+^{k}\big]$ holds. Finally, for
every $c>0$ we will let $\gamma^{(c)}\in\cP$ denote the distribution of
$T/c$, i.e., the push-forward of $\gamma$ via $t\mapsto t/c$.

We start with the following simple lemma describing the effect of rescaling
a Williamson measure.

\begin{lemma}\label{lem:scaling}
Suppose that $\gamma\in\cP$ and $c>0$. Then $G_{\gamma^{(c)}}(s)=G_\gamma(cs)$
and $\Pk{k}_{\gamma^{(c)}}(s)=\Pk{k}_\gamma(cs)$ hold for all $s>0$ and
$k\in\N$. Moreover, we have $s_{\gamma^{(c)}}=s_\gamma/c$ as well as
$C_{\gamma^{(c)}}=C_\gamma$.
\end{lemma}

\begin{proof}
Let $T\sim\gamma$. Then $G_{\gamma^{(c)}}(s)=\PP(T/c\leq s)=G_\gamma(cs)$ and
$\Pk{k}_{\gamma^{(c)}}(s)=\EE[(1-T/(cs))_+^k]=\Pk{k}_\gamma(cs)$, from which
$s_{\gamma^{(c)}}=s_\gamma/c$ follows immediately. In particular we have
$\psi_{\gamma^{(c)}}(z)=\psi_\gamma(z/c)$ for every $z\geq 0$, so the
pseudo-inverse of $\psi_{\gamma^{(c)}}$ is given by $c\,\varphi$, whereby
$\varphi$ denotes the pseudo-inverse of $\psi_\gamma$. Having this,
$C_{\gamma^{(c)}}(\bx)=\psi_\gamma\big(\tfrac1c\sum_{i=1}^d c\,\varphi(x_i)\big)
=C_\gamma(\bx)$ for every $\bx\in\II^d$ follows, which completes the proof.
\end{proof}

The next lemma gathers some properties of the transforms $\Pk{k}_\gamma$
which will be used frequently in the sequel.

\begin{lemma}\label{lem:Psi}
Suppose that $\gamma\in\cP$, let $G:=G_\gamma$ and $k\in\N$. Then the
following assertions hold:
\begin{enumerate}
\item\label{it:Psi1} For every $s>0$ we have
\begin{equation}\label{eq:PsiG}
  \Pk{k}_\gamma(s)=\frac{k}{s^{k}}\int_0^s (s-t)^{k-1}\,G(t)\,\dd t .
\end{equation}
\item\label{it:Psi2} $\Pk{k}_\gamma$ is continuous and non-decreasing on
$(0,\infty)$ with values in $[0,1)$. Moreover, $\Pk{k}_\gamma(s)=0$ holds if,
and only if, $s\leq s_\gamma$, the function $\Pk{k}_\gamma$ is strictly
increasing on $[s_\gamma,\infty)\cap(0,\infty)$, and we have
$\lim_{s\to\infty}\Pk{k}_\gamma(s)=1$ as well as
$\lim_{s\downarrow 0}\Pk{k}_\gamma(s)=0$.
\end{enumerate}
\end{lemma}

\begin{proof}
\ref{it:Psi1} Considering that for $0<t<s$ we have
$(1-t/s)^k=\int_t^s \frac{k}{s}(1-u/s)^{k-1}\,\dd u$, applying Fubini's
theorem yields
\[
  \Pk{k}_\gamma(s)=\int_0^s\frac{k}{s}\Big(1-\frac{u}{s}\Big)^{k-1}
  \gamma((0,u))\,\dd u ,
\]
from which eq.~\eqref{eq:PsiG} follows since $\gamma((0,u))=G(u^-)$ coincides
with $G(u)$ for all but at most countably many $u$.

\ref{it:Psi2} Continuity is a direct consequence of eq.~\eqref{eq:PsiG}. For
$0<s<s'$ and $t>0$ we have $(1-t/s)_+\leq(1-t/s')_+$ with strict inequality
whenever $t<s'$. Hence $\Pk{k}_\gamma$ is non-decreasing and
$\Pk{k}_\gamma(s)<\Pk{k}_\gamma(s')$ holds whenever $\gamma((0,s'))>0$, i.e.,
whenever $s'>s_\gamma$. Furthermore $\Pk{k}_\gamma(s)=0$ is equivalent to
$\gamma((0,s))=0$, i.e., to $s\leq s_\gamma$, and $\Pk{k}_\gamma(s)<1$ follows
from $(1-t/s)_+^k<1$ for $t>0$. Finally, dominated convergence yields
$\lim_{s\to\infty}\Pk{k}_\gamma(s)=1$, and the inequality
$\Pk{k}_\gamma(s)\leq\gamma((0,s))$ implies
$\lim_{s\downarrow0}\Pk{k}_\gamma(s)=0$.
\end{proof}

We now turn to the Kendall distribution function and express it in terms of
the Williamson measure (see \cite[Theorem~3.1 and Section~4.2]{MN}; notice
that for $d\geq 3$ the result is equivalent to \cite[Theorem~5.5]{KDT}).

\begin{lemma}\label{lem:kendall}
Suppose that $\gamma\in\cP$, let $T$ be a random variable with $T\sim\gamma$
and let $\FK$ denote the Kendall distribution function of $C_\gamma$. Then
$\FK$ is the distribution function of the random variable
$\Psi_\gamma(T)=\psi_\gamma(1/T)$. Moreover, $C_\gamma$ is strict if, and only
if, $s_\gamma=0$.
\end{lemma}

\begin{proof}
Let $\varphi$ denote the pseudo-inverse of $\psi:=\psi_\gamma$, let
$\bU=(U_1,\ldots,U_d)\sim C_\gamma$ and set $R:=1/T$. Then
$\psi(z)=\EE[(1-z/R)_+^{\,n}]$ holds for every $z\geq 0$, i.e., $\psi$ is the
Williamson $d$-transform of the distribution of $R$ in the sense of \cite{MN}.
According to \cite[Theorem~3.1]{MN} the random vector
$(\varphi(U_1),\ldots,\varphi(U_d))$ has the same distribution as
$R\,\mathbf{S}$, whereby $\mathbf{S}$ is uniformly distributed on the unit
simplex $\{\bx\in[0,\infty)^d\colon x_1+\cdots+x_d=1\}$ and independent of
$R$. Consequently, $\varphi(U_1)+\cdots+\varphi(U_d)$ has the same
distribution as $R$, from which it follows that
$C_\gamma(\bU)=\psi(\varphi(U_1)+\cdots+\varphi(U_d))$ has the same
distribution as $\psi(R)=\Psi_\gamma(T)$. Finally, considering
$\psi_\gamma(z)=\Psi_\gamma(1/z)$ and applying \cref{lem:Psi} yields
$\varphi(0)=\inf\{z\in[0,\infty]\colon\psi_\gamma(z)=0\}=1/s_\gamma$ (with
the convention $1/0:=\infty$), which proves the second assertion.
\end{proof}

For our purposes the following reformulation of \cref{lem:kendall} will be
key.

\begin{corollary}\label{cor:key}
Suppose that $\gamma\in\cP$ and let $\FK$ denote the Kendall distribution
function of $C_\gamma$. Then the identity
\begin{equation}\label{eq:key}
  G_\gamma(s)=\FK\big(\Psi_\gamma(s)\big)
\end{equation}
holds for every $s>0$ fulfilling $s\geq s_\gamma$. In particular,
eq.~\eqref{eq:key} holds for every $s>0$ if $C_\gamma$ is strict.
\end{corollary}

\begin{proof}
Let $T\sim\gamma$ and fix $s>0$ with $s\geq s_\gamma$. Considering that
$T\geq s_\gamma$ and $T>0$ hold almost surely and that, according to
\cref{lem:Psi}, $\Psi_\gamma$ is strictly increasing on
$[s_\gamma,\infty)\cap(0,\infty)$, the events
$\{\Psi_\gamma(T)\leq\Psi_\gamma(s)\}$ and $\{T\leq s\}$ coincide up to a null
set. Applying \cref{lem:kendall} therefore yields
$\FK(\Psi_\gamma(s))=\PP(\Psi_\gamma(T)\leq\Psi_\gamma(s))=\PP(T\leq s)
=G_\gamma(s)$, which completes the proof.
\end{proof}

\begin{remark}\label{rem:classical}
Notice that eq.~\eqref{eq:key} can also be verified directly via the
well-known formula for the Kendall distribution function in terms of the
derivatives of the generator (see \cite[Proposition~4.5]{MN} as well as
\cite[eq.~(3.7)]{KDT}), according to which
\[
  \FK\big(\psi(z)\big)=\sum_{k=0}^{n-1}\frac{(-z)^{k}}{k!}\,\psi^{(k)}(z)
  +\frac{(-z)^{n}}{n!}\,D^-\psi^{(n-1)}(z)
\]
holds for $\psi=\psi_\gamma$ and every $z\in(0,\varphi(0))$, whereby
$D^-\psi^{(n-1)}$ denotes the left-hand derivative of $\psi^{(n-1)}$. In
fact, letting $T\sim\gamma$ and differentiating eq.~\eqref{eq:williamson}
yields
\[
  \psi^{(k)}(z)=(-1)^{k}\,\frac{n!}{(n-k)!}\,\EE\big[T^{k}(1-zT)_+^{n-k}\big],
  \qquad k\in\{0,\ldots,n-1\},
\]
as well as
$D^-\psi^{(n-1)}(z)=(-1)^{n}\,n!\,\EE\big[T^{n}\,\ind_{(0,1]}(zT)\big]$ (for
the latter identity also see \cite[Lemma~5.3]{KDT}), so applying the binomial
theorem we get
\[
  \FK\big(\psi(z)\big)=\EE\Big[\ind_{(0,1]}(zT)\sum_{k=0}^{n}\binom{n}{k}
  (zT)^{k}(1-zT)^{n-k}\Big]=\PP\big(T\leq\tfrac{1}{z}\big)
  =G_\gamma\big(\tfrac{1}{z}\big),
\]
which is eq.~\eqref{eq:key} for $s=1/z>s_\gamma$.
\end{remark}

\section{A cyclic pursuit system}\label{sec:pursuit}

Throughout this section and \cref{sec:proof} we consider two measures
$\gamma_1,\gamma_2\in\cP$ fulfilling that the corresponding Archimedean
copulas $C_{\gamma_1}$ and $C_{\gamma_2}$ have the same Kendall distribution
function, which will be denoted by $K$. In order to keep notation as simple
as possible, for $i\in\{1,2\}$ and $k\in\N$ we will write
\[
  G_i:=G_{\gamma_i},\qquad s_i:=s_{\gamma_i},\qquad
  \Pk{k}_i:=\Pk{k}_{\gamma_i},\qquad \Pk{0}_i:=G_i .
\]
We start with the following lemma providing a differentiation rule for the
transforms $\Pk{k}_\gamma$, which will be key in the sequel.

\begin{lemma}\label{lem:diff}
Suppose that $\gamma\in\cP$ and $k\in\N$, and set $\Pk{0}_\gamma:=G_\gamma$.
Then the following assertions hold:
\begin{enumerate}
\item\label{it:diff1} The restriction of $\Pk{k}_\gamma$ to
$(s_\gamma,\infty)$ is an increasing homeomorphism of $(s_\gamma,\infty)$ onto
$(0,1)$.
\item\label{it:diff2} $\Pk{k}_\gamma$ has a right-hand derivative at every
$s>0$, which is given by
\begin{equation}\label{eq:diff}
  \big(\Pk{k}_\gamma\big)'_+(s)
  =\frac{k}{s}\Big(\Pk{k-1}_\gamma(s)-\Pk{k}_\gamma(s)\Big)
\end{equation}
and which is strictly positive for every $s>s_\gamma$. For $k\geq 2$ the
function $\Pk{k}_\gamma$ is even continuously differentiable on $(0,\infty)$.
\end{enumerate}
\end{lemma}

\begin{proof}
\ref{it:diff1} According to \cref{lem:Psi} the function $\Pk{k}_\gamma$ is
continuous on $(0,\infty)$ and strictly increasing on
$[s_\gamma,\infty)\cap(0,\infty)$, and it fulfills
$\lim_{s\to\infty}\Pk{k}_\gamma(s)=1$ as well as
$\lim_{s\downarrow s_\gamma}\Pk{k}_\gamma(s)=0$, from which the first
assertion follows immediately.

\ref{it:diff2} First notice that for $k\geq 2$ the function
$s\mapsto\int_0^s(s-t)^{k-1}G_\gamma(t)\,\dd t$ is continuously
differentiable with derivative $(k-1)\int_0^s(s-t)^{k-2}G_\gamma(t)\,\dd t$,
and that for $k=1$ the function $s\mapsto\int_0^sG_\gamma(t)\,\dd t$ has
right-hand derivative $G_\gamma(s)$ at every $s>0$ since $G_\gamma$ is
right-continuous. Having this, eq.~\eqref{eq:diff} follows from
eq.~\eqref{eq:PsiG} and the product rule, and for $k\geq 2$ the right-hand
side of eq.~\eqref{eq:diff} is continuous. Finally,
\[
  \Pk{k-1}_\gamma(s)-\Pk{k}_\gamma(s)
  =\int_{(0,s)}\Big(1-\frac{t}{s}\Big)^{k-1}\frac{t}{s}\,\dd\gamma(t)
  +\ind_{\{1\}}(k)\,\gamma(\{s\})
\]
is strictly positive for every $s>s_\gamma$ since in this case we have
$\gamma((0,s))>0$.
\end{proof}

According to \cref{lem:diff}, for every $k\in\{1,\ldots,n\}$ the function
\[
  \sigma_k:=\big(\Pk{k}_2|_{(s_2,\infty)}\big)^{-1}\circ
  \Pk{k}_1|_{(s_1,\infty)}
\]
is an increasing homeomorphism of $(s_1,\infty)$ onto $(s_2,\infty)$. In other
words: for every $s>s_1$ the value $\sigma_k(s)$ is the unique point $t>0$
fulfilling $\Pk{k}_2(t)=\Pk{k}_1(s)$. Furthermore we set
\[
  \varrho_k(s):=\frac{\sigma_k(s)}{s},\qquad s>s_1,
  \qquad\text{as well as}\qquad \sigma_0:=\sigma_n,\quad
  \varrho_0:=\varrho_n .
\]
Considering \cref{lem:scaling}, $\varrho_k(s)$ is the unique scaling factor
$c>0$ fulfilling $\Pk{k}_{\gamma_1}(s)=\Pk{k}_{\gamma_2^{(c)}}(s)$, i.e., the
factor by which $\gamma_2$ has to be rescaled in order to match the transform
of order $k$ of $\gamma_1$ at the point $s$. In particular,
$\gamma_1=\gamma_2^{(c)}$ implies $\varrho_k(s)=c$ for all $k\in\{1,\ldots,n\}$
and $s>s_1$. The assumption that the Kendall distribution functions of
$C_{\gamma_1}$ and $C_{\gamma_2}$ coincide enters via the following simple
identity, which, loosely speaking, closes the cycle.

\begin{lemma}\label{lem:cycle}
The identity $G_1=G_2\circ\sigma_n$ holds on $(s_1,\infty)$. As a consequence,
$\Pk{k}_1=\Pk{k}_2\circ\sigma_k$ holds on $(s_1,\infty)$ for every
$k\in\{0,1,\ldots,n\}$.
\end{lemma}

\begin{proof}
Fix $s>s_1$. Considering $\sigma_n(s)>s_2$, applying \cref{cor:key} to
$\gamma_1$ as well as to $\gamma_2$ and using the definition of $\sigma_n$
yields
\[
  G_1(s)=K\big(\Pk{n}_1(s)\big)=K\big(\Pk{n}_2(\sigma_n(s))\big)
  =G_2(\sigma_n(s)).
\]
For $k\geq 1$ the second assertion holds by definition of $\sigma_k$, and for
$k=0$ it coincides with the first one.
\end{proof}

Building upon \cref{lem:diff} and \cref{lem:cycle} we can now derive the
main result of this section, saying that each of the functions $\varrho_k$
pursues its predecessor.

\begin{proposition}[Cyclic pursuit]\label{prop:pursuit}
For every $k\in\{1,\ldots,n\}$ the function $\varrho_k$ is continuous and has
a right-hand derivative at every $s>s_1$ $($for $k\geq 2$ it is even
continuously differentiable$)$. Moreover, the identity
\begin{equation}\label{eq:pursuit}
  \frac{s\,(\varrho_k)'_+(s)}{\varrho_k(s)}
  =\frac{\Pk{k-1}_2\big(\sigma_{k-1}(s)\big)-\Pk{k-1}_2\big(\sigma_k(s)\big)}
        {\Pk{k-1}_2\big(\sigma_k(s)\big)-\Pk{k}_2\big(\sigma_k(s)\big)}
\end{equation}
holds, whereby the denominator is strictly positive. In particular, for all
$k\in\{1,\ldots,n\}$ and $s>s_1$ the following implications hold:
\begin{align*}
  \varrho_{k-1}(s)\leq\varrho_k(s)\ &\Longrightarrow\ (\varrho_k)'_+(s)\leq 0,\\
  \varrho_{k-1}(s)\geq\varrho_k(s)\ &\Longrightarrow\ (\varrho_k)'_+(s)\geq 0.
\end{align*}
\end{proposition}

\begin{proof}
Fix $k\in\{1,\ldots,n\}$ and $s>s_1$, and recall that $\sigma_k(s)>s_2$.
According to \cref{lem:diff} the increasing homeomorphism
$\Pk{k}_2|_{(s_2,\infty)}$ has a strictly positive right-hand derivative at
$\sigma_k(s)$, so its inverse has right-hand derivative
$1/(\Pk{k}_2)'_+(\sigma_k(s))$ at the point
$\Pk{k}_2(\sigma_k(s))=\Pk{k}_1(s)$. Considering that $\Pk{k}_1$ is strictly
increasing on $(s_1,\infty)$ and has a right-hand derivative at $s$, the
chain rule for right-hand derivatives and eq.~\eqref{eq:diff} yield
\[
  (\sigma_k)'_+(s)=\frac{(\Pk{k}_1)'_+(s)}{(\Pk{k}_2)'_+(\sigma_k(s))}
  =\frac{\sigma_k(s)}{s}\cdot
  \frac{\Pk{k-1}_1(s)-\Pk{k}_1(s)}
       {\Pk{k-1}_2(\sigma_k(s))-\Pk{k}_2(\sigma_k(s))} .
\]
According to \cref{lem:cycle} we have
$\Pk{k-1}_1(s)=\Pk{k-1}_2(\sigma_{k-1}(s))$ as well as
$\Pk{k}_1(s)=\Pk{k}_2(\sigma_k(s))$. Having this, eq.~\eqref{eq:pursuit}
follows from
\[
  \frac{s\,(\varrho_k)'_+(s)}{\varrho_k(s)}
  =\frac{s\,(\sigma_k)'_+(s)}{\sigma_k(s)}-1 ,
\]
and positivity of the denominator has already been established in
\cref{lem:diff}. For $k\geq 2$ all involved functions are continuously
differentiable. Finally, the two implications are a direct consequence of
the fact that $\Pk{k-1}_2$ is non-decreasing (recall $\Pk{0}_2=G_2$ and
$\sigma_0=\sigma_n$).
\end{proof}

Loosely speaking, \cref{prop:pursuit} states that each $\varrho_k$ moves
towards its predecessor $\varrho_{k-1}$, whereby $\varrho_1$ moves towards
$\varrho_n$. As a direct consequence, the system contracts as $s$ increases:

\begin{corollary}\label{cor:monotone}
Define $M,m\colon(s_1,\infty)\to(0,\infty)$ by
\[
  M(s):=\max_{k\in\{1,\ldots,n\}}\varrho_k(s),\qquad
  m(s):=\min_{k\in\{1,\ldots,n\}}\varrho_k(s).
\]
Then the following assertions hold:
\begin{enumerate}
\item\label{it:mon1} $M$ is non-increasing and $m$ is non-decreasing.
\item\label{it:mon2} If $\sup_{s>s_1}M(s)\leq\inf_{s>s_1}m(s)$, then
$c:=\sup_{s>s_1}M(s)\in(0,\infty)$ and $\gamma_1=\gamma_2^{(c)}$.
\end{enumerate}
\end{corollary}

\begin{proof}
\ref{it:mon1} We will use the following elementary fact: if
$f\colon(a,b)\to\R$ is continuous and has a right-hand derivative fulfilling
$f'_+\leq 0$ at every point, then $f$ is non-increasing. (Suppose, on the
contrary, that $f(y)>f(x)$ for some $x<y$, choose $\epsilon>0$ such that
$g(z):=f(z)-f(x)-\epsilon(z-x)$ fulfills $g(y)>0$, and let $z^*$ denote the
largest point in $[x,y]$ with $g(z^*)\leq 0$. Then $z^*<y$ and $g>g(z^*)$ on
$(z^*,y]$, hence $g'_+(z^*)\geq 0$, which contradicts
$g'_+(z^*)=f'_+(z^*)-\epsilon<0$.)

The functions $\varrho_1,\ldots,\varrho_n$ are continuous and have right-hand
derivatives at every $s>s_1$. Hence $M$ is continuous and has right-hand
derivative
\[
  M'_+(s)=\max\big\{(\varrho_k)'_+(s)\colon k\in\{1,\ldots,n\},\
  \varrho_k(s)=M(s)\big\}.
\]
If $\varrho_k(s)=M(s)$, then $\varrho_{k-1}(s)\leq\varrho_k(s)$ (recall
$\varrho_0=\varrho_n$), so \cref{prop:pursuit} yields
$(\varrho_k)'_+(s)\leq 0$. Consequently, we have $M'_+\leq 0$ on
$(s_1,\infty)$, implying that $M$ is non-increasing. The proof for $m$ is
analogous (consider $-m$).

\ref{it:mon2} Considering $m\leq M$ the assumption implies $M(s)=m(s)=c$ for
every $s>s_1$, so $c\in(0,\infty)$ and $\sigma_k(s)=c\,s$ for all $s>s_1$
and $k\in\{1,\ldots,n\}$. As $\sigma_n$ maps $(s_1,\infty)$ onto
$(s_2,\infty)$ we get $s_2=c\,s_1$, and \cref{lem:cycle} yields
$G_1(s)=G_2(cs)$ for every $s>s_1$. If $s_1>0$, then for $s\in(0,s_1)$ we
have $G_1(s)=0=G_2(cs)$, and $G_1(s_1)=G_2(cs_1)$ follows from
right-continuity. Altogether $G_1(s)=G_2(cs)$ holds for every $s>0$, which,
according to \cref{lem:scaling}, is equivalent to $\gamma_1=\gamma_2^{(c)}$.
\end{proof}

\section{Proof of the main result}\label{sec:proof}

We keep the assumptions and the notation from \cref{sec:pursuit}, i.e.,
$\gamma_1,\gamma_2\in\cP$ induce the same Kendall distribution function $K$,
and tackle the non-strict and the strict case separately.

\subsection{The non-strict case}\label{subsec:nonstrict}

In the non-strict case the monotonicity properties established in
\cref{cor:monotone} directly yield the desired result:

\begin{proposition}\label{prop:nonstrict}
If $s_1>0$ or $s_2>0$, then both $s_1$ and $s_2$ are strictly positive and
$\gamma_1=\gamma_2^{(c)}$ holds for $c:=s_2/s_1$.
\end{proposition}

\begin{proof}
Suppose first that $s_1>0$. Being an increasing homeomorphism of
$(s_1,\infty)$ onto $(s_2,\infty)$, each $\sigma_k$ fulfills
$\lim_{s\downarrow s_1}\sigma_k(s)=s_2$, hence
$\lim_{s\downarrow s_1}\varrho_k(s)=s_2/s_1$ for every $k\in\{1,\ldots,n\}$
and therefore
\[
  \lim_{s\downarrow s_1}M(s)=\frac{s_2}{s_1}=\lim_{s\downarrow s_1}m(s).
\]
Since $M$ is non-increasing and $m$ is non-decreasing this implies
$\sup_{s>s_1}M(s)=s_2/s_1=\inf_{s>s_1}m(s)$, and applying \cref{cor:monotone}
yields $s_2/s_1\in(0,\infty)$ as well as $\gamma_1=\gamma_2^{(c)}$ for
$c=s_2/s_1$. If $s_2>0$, then applying what has already been proved to the
pair $(\gamma_2,\gamma_1)$ yields $s_1/s_2\in(0,\infty)$, i.e., $s_1>0$, and
we are in the situation considered before.
\end{proof}

\subsection{The strict case: a priori bounds and pursuit rates}
\label{subsec:bounds}

In the remainder of this section we assume $s_1=s_2=0$, i.e., that both
$C_{\gamma_1}$ and $C_{\gamma_2}$ are strict (see \cref{lem:kendall}). Then
$\sigma_1,\ldots,\sigma_n$ are increasing homeomorphisms of $(0,\infty)$, the
functions $\varrho_1,\ldots,\varrho_n$, $M$ and $m$ are defined on
$(0,\infty)$, and we set
\[
  M^*:=\sup_{s>0}M(s)\in(0,\infty],\qquad
  m_*:=\inf_{s>0}m(s)\in[0,\infty).
\]
Contrary to the non-strict case the pursuit system now has no starting
point, so additional arguments are needed. We start with the following simple
lemma relating transforms of consecutive order.

\begin{lemma}\label{lem:sandwich}
Suppose that $\gamma\in\cP$ and let $k\geq 2$ be an integer. Then
\[
  \Pk{k}_\gamma(s)\ \leq\ \Pk{k-1}_\gamma(s)\ \leq\
  \Pk{k}_\gamma\big(\tfrac{k}{k-1}\,s\big)
\]
holds for every $s>0$.
\end{lemma}

\begin{proof}
Let $T\sim\gamma$. The first inequality follows from $(1-T/s)_+\in[0,1]$.
Applying the inequality between the geometric and the arithmetic mean to the
$k$ numbers $1-u,\ldots,1-u,1$ (with $1-u$ appearing $k-1$ times) yields
$(1-u)^{k-1}\leq\big(1-\tfrac{k-1}{k}\,u\big)^{k}$ for every $u\in[0,1]$, and
therefore
$(1-u)_+^{k-1}\leq\big(1-\tfrac{k-1}{k}\,u\big)_+^{k}$ for every $u\geq 0$.
Considering $u=T/s$ and taking expectations we get
\[
  \Pk{k-1}_\gamma(s)\leq
  \EE\Big[\Big(1-\frac{k-1}{k}\cdot\frac{T}{s}\Big)_+^{k}\Big]
  =\Pk{k}_\gamma\big(\tfrac{k}{k-1}\,s\big),
\]
which completes the proof.
\end{proof}

Having this, we can derive a priori bounds for the functions $\sigma_k$ and
$\varrho_k$.

\begin{lemma}[A priori bounds]\label{lem:apriori}
The following assertions hold:
\begin{enumerate}
\item\label{it:ap1} For every $k\in\{2,\ldots,n\}$ and every $s>0$ we have
\[
  \sigma_k(s)\leq\tfrac{k}{k-1}\,\sigma_{k-1}(s)
  \qquad\text{and}\qquad
  \sigma_{k-1}(s)\leq\sigma_k\big(\tfrac{k}{k-1}\,s\big).
\]
\item\label{it:ap2} $M(s)\leq n^2\,m(ns)$ for every $s>0$.
\item\label{it:ap3} $M^*=\lim_{s\downarrow 0}M(s)$,
$m_*=\lim_{s\downarrow 0}m(s)$ and $0<m_*\leq M^*\leq n^2\,m_*<\infty$.
\end{enumerate}
\end{lemma}

\begin{proof}
\ref{it:ap1} Applying \cref{lem:cycle} and \cref{lem:sandwich} (to $\gamma_1$
as well as to $\gamma_2$) yields
\[
  \Pk{k}_2\big(\sigma_k(s)\big)=\Pk{k}_1(s)\leq\Pk{k-1}_1(s)
  =\Pk{k-1}_2\big(\sigma_{k-1}(s)\big)
  \leq\Pk{k}_2\big(\tfrac{k}{k-1}\,\sigma_{k-1}(s)\big)
\]
as well as
\[
  \Pk{k-1}_2\big(\sigma_{k-1}(s)\big)=\Pk{k-1}_1(s)
  \leq\Pk{k}_1\big(\tfrac{k}{k-1}\,s\big)
  =\Pk{k}_2\big(\sigma_k(\tfrac{k}{k-1}\,s)\big)
  \leq\Pk{k-1}_2\big(\sigma_k(\tfrac{k}{k-1}\,s)\big).
\]
Considering that $\Pk{k}_2$ and $\Pk{k-1}_2$ are strictly increasing on
$(0,\infty)$, both inequalities follow.

\ref{it:ap2} Let $j,k\in\{1,\ldots,n\}$ and $s>0$. Iterating the first
inequality in \ref{it:ap1} yields $\sigma_n(t)\leq\frac{n}{j}\,\sigma_j(t)$
for every $t>0$, and iterating the second one yields
$\sigma_k(s)\leq\sigma_n\big(\frac{n}{k}\,s\big)$. Consequently,
\[
  \varrho_k(s)=\frac{\sigma_k(s)}{s}
  \leq\frac{\sigma_n\big(\frac{n}{k}\,s\big)}{s}
  \leq\frac{n}{j}\cdot\frac{\sigma_j\big(\frac{n}{k}\,s\big)}{s}
  =\frac{n^2}{jk}\,\varrho_j\big(\tfrac{n}{k}\,s\big)
  \leq n^2\,\varrho_j\big(\tfrac{n}{k}\,s\big).
\]
Since $j$ was arbitrary and since $m$ is non-decreasing according to
\cref{cor:monotone}, we get
$\varrho_k(s)\leq n^2\,m\big(\frac{n}{k}\,s\big)\leq n^2\,m(ns)$ for every
$k\in\{1,\ldots,n\}$, i.e., $M(s)\leq n^2\,m(ns)$.

\ref{it:ap3} The first two identities are direct consequences of
\cref{cor:monotone}. Considering $s\downarrow 0$ in \ref{it:ap2} yields
$M^*\leq n^2\,m_*$. Since $M^*\geq M(1)>0$ and $m_*\leq m(1)<\infty$ this
shows $m_*>0$ as well as $M^*<\infty$, and $m_*\leq M^*$ follows from
$m\leq M$.
\end{proof}

The next lemma shows that for $k\geq 2$ the factor $\varrho_k$ approaches its
predecessor $\varrho_{k-1}$ from above at a rate which is bounded from below.

\begin{lemma}[Pursuit rate]\label{lem:rate}
Let $k\in\{2,\ldots,n\}$ and suppose that $s>0$ fulfills
$\varrho_k(s)>\varrho_{k-1}(s)$. Then
\[
  \frac{s\,\varrho_k'(s)}{\varrho_k(s)}\ \leq\
  -\Big(\frac{\varrho_k(s)}{\varrho_{k-1}(s)}-1\Big)\ \leq\
  -\ln\frac{\varrho_k(s)}{\varrho_{k-1}(s)} .
\]
\end{lemma}

\begin{proof}
Set $t:=\sigma_k(s)$ as well as
$\lambda:=\varrho_k(s)/\varrho_{k-1}(s)=\sigma_k(s)/\sigma_{k-1}(s)$, so that
$\sigma_{k-1}(s)=t/\lambda$ and, according to \cref{lem:apriori},
$1<\lambda\leq\frac{k}{k-1}\leq 2$. Letting $T\sim\gamma_2$ and $U:=T/t$ we
have $\Pk{j}_2(t)=\EE[(1-U)_+^{j}]$ and
$\Pk{j}_2(t/\lambda)=\EE[(1-\lambda U)_+^{j}]$ for every $j\in\N$.
Considering $k-1\geq 1$ and using the fact that
$(1-U)_+^{k-1}-(1-U)_+^{k}=U\,(1-U)_+^{k-1}$ holds, eq.~\eqref{eq:pursuit}
therefore reads
\begin{equation}\label{eq:rate}
  \frac{s\,\varrho_k'(s)}{\varrho_k(s)}
  =-\,\frac{\EE\big[(1-U)_+^{k-1}-(1-\lambda U)_+^{k-1}\big]}
           {\EE\big[U\,(1-U)_+^{k-1}\big]} .
\end{equation}
We now show that
\begin{equation}\label{eq:pointwise}
  (1-u)_+^{k-1}-(1-\lambda u)_+^{k-1}\ \geq\ (\lambda-1)\,u\,(1-u)_+^{k-1}
\end{equation}
holds for every $u\geq 0$. For $u\geq 1$ both sides vanish. In case of
$\frac{1}{\lambda}\leq u<1$ the left-hand side equals $(1-u)^{k-1}$ and
we have $(\lambda-1)\,u\leq\lambda-1\leq 1$. Finally, for
$u<\frac{1}{\lambda}$ the numbers $a:=1-u$ and $b:=1-\lambda u$ fulfill
$0<b\leq a\leq 1$ as well as $a-b=(\lambda-1)\,u$, and therefore
\[
  a^{k-1}-b^{k-1}=(a-b)\sum_{j=0}^{k-2}a^{j}\,b^{k-2-j}
  \ \geq\ (a-b)\,a^{k-2}\ \geq\ (a-b)\,a^{k-1},
\]
which proves eq.~\eqref{eq:pointwise}. Considering $u=U$ in
eq.~\eqref{eq:pointwise}, taking expectations and using that the denominator
in eq.~\eqref{eq:rate} is strictly positive yields the first inequality, and
the second one follows from $\ln\lambda\leq\lambda-1$.
\end{proof}

\subsection{The strict case: a cascade argument}\label{subsec:cascade}

It is convenient to pass to logarithmic coordinates. According to
\cref{lem:apriori} we have $M^*\in(0,\infty)$, and for $k\in\{1,\ldots,n\}$
and $x\in\R$ we set
\[
  \tau_k(x):=\ln\varrho_k(e^{x}),\qquad \ell:=\ln M^*,\qquad
  p_k(x):=\ell-\tau_k(x).
\]
The subsequent lemma translates the results established so far to this
setting.

\begin{lemma}\label{lem:log}
The functions $p_1,\ldots,p_n\colon\R\to[0,\infty)$ are continuous and have
the following properties:
\begin{enumerate}
\item\label{it:log1} $p_k(x)\leq p_k(y)+(x-y)$ holds for every
$k\in\{1,\ldots,n\}$ and all $y\leq x$.
\item\label{it:log2}
$\displaystyle\lim_{x\to-\infty}\,\min_{k\in\{1,\ldots,n\}}p_k(x)=0$.
\item\label{it:log3} For every $k\in\{2,\ldots,n\}$ the function $p_k$ is
continuously differentiable, and
\[
  p_k'(x)\geq p_{k-1}(x)-p_k(x)
\]
holds for every $x\in\R$ fulfilling $p_{k-1}(x)>p_k(x)$.
\end{enumerate}
\end{lemma}

\begin{proof}
Continuity is clear and $p_k\geq 0$ follows from $\varrho_k\leq M\leq M^*$.
Since the function $x\mapsto x+\tau_k(x)=\ln\sigma_k(e^{x})$ is increasing,
we have $\tau_k(x)\geq\tau_k(y)-(x-y)$ for all $y\leq x$, which is
\ref{it:log1}. The second assertion follows from
$\min_kp_k(x)=\ell-\ln M(e^{x})$ and \cref{lem:apriori}. Finally, let
$k\geq 2$. According to \cref{prop:pursuit} the function $\tau_k$ is
continuously differentiable with $\tau_k'(x)=s\,\varrho_k'(s)/\varrho_k(s)$
for $s=e^{x}$. If $p_{k-1}(x)>p_k(x)$, then $\varrho_k(s)>\varrho_{k-1}(s)$
and \cref{lem:rate} yields
\[
  p_k'(x)=-\tau_k'(x)\geq\ln\frac{\varrho_k(s)}{\varrho_{k-1}(s)}
  =\tau_k(x)-\tau_{k-1}(x)=p_{k-1}(x)-p_k(x),
\]
which completes the proof.
\end{proof}

The key step is the following lemma: if $\varrho_k$ is close to the supremum
$M^*$ at some point, then so is its predecessor $\varrho_{k-1}$. In fact,
otherwise $\varrho_k$ would have been pulled away from $M^*$ shortly before.

\begin{lemma}[Cascade]\label{lem:cascade}
Let $k\in\{2,\ldots,n\}$ and $x\in\R$. Then
\begin{equation}\label{eq:cascade}
  p_{k-1}(x)\leq\frac{p_k(x)}{1-e^{-h}}+h
\end{equation}
holds for every $h>0$. In particular, $p_k(x)\leq 1$ implies
$p_{k-1}(x)\leq 3\sqrt{p_k(x)}$.
\end{lemma}

\begin{proof}
Fix $h>0$ and set $c:=\min\{p_{k-1}(y)\colon y\in[x-h,x]\}\geq 0$. According
to \cref{lem:log} we have $p_{k-1}(x)\leq p_{k-1}(y)+h$ for every
$y\in[x-h,x]$, hence $p_{k-1}(x)\leq c+h$, and it suffices to show
\begin{equation}\label{eq:cascade2}
  p_k(x)\geq\big(1-e^{-h}\big)\,c .
\end{equation}
Consider the continuously differentiable function $q:=c-p_k$ on $[x-h,x]$.
If $q(x)\leq 0$, then $p_k(x)\geq c$ and eq.~\eqref{eq:cascade2} holds.
Suppose now that $q(x)>0$ and set
\[
  a:=\inf\big\{y\in[x-h,x]\colon q>0\text{ on }[y,x]\big\},
\]
so that $q>0$ on $(a,x]$. For every $y\in(a,x]$ we have
$p_k(y)<c\leq p_{k-1}(y)$, so \cref{lem:log} yields
\[
  q'(y)=-p_k'(y)\leq-\big(p_{k-1}(y)-p_k(y)\big)\leq-\big(c-p_k(y)\big)
  =-q(y).
\]
Consequently, the function $y\mapsto e^{y}q(y)$ is non-increasing on $[a,x]$,
which implies $e^{a}q(a)\geq e^{x}q(x)>0$ and therefore $q(a)>0$. If
$a>x-h$ held, then continuity of $q$ would yield $q>0$ on $[a-\epsilon,x]$
for some $\epsilon\in(0,a-(x-h)]$, contradicting the definition of $a$.
Hence $a=x-h$, and using $p_k(x-h)\geq 0$ we obtain
\[
  q(x)\leq e^{-h}\,q(x-h)\leq e^{-h}\,c ,
\]
i.e., $p_k(x)=c-q(x)\geq(1-e^{-h})\,c$, which proves
eq.~\eqref{eq:cascade2}.

To prove the second assertion suppose that $p_k(x)\leq 1$. If $p_k(x)=0$,
then considering $h\downarrow 0$ in eq.~\eqref{eq:cascade} yields
$p_{k-1}(x)=0$. If $p_k(x)>0$, then we choose $h:=\sqrt{p_k(x)}\in(0,1]$.
Since $1-e^{-h}\geq h-\frac{h^2}{2}\geq\frac{h}{2}$, eq.~\eqref{eq:cascade}
yields $p_{k-1}(x)\leq 2\sqrt{p_k(x)}+\sqrt{p_k(x)}=3\sqrt{p_k(x)}$.
\end{proof}

Applying \cref{lem:cascade} iteratively yields the following result, which
is the main ingredient for the proof of \cref{thm:strict}.

\begin{proposition}\label{prop:top}
We have $\displaystyle\lim_{s\downarrow 0}\varrho_1(s)=M^*$.
\end{proposition}

\begin{proof}
We have to show $\lim_{x\to-\infty}p_1(x)=0$. For $n=1$ this is the second
assertion of \cref{lem:log}, so suppose that $n\geq 2$. Fix
$\varepsilon\in(0,1]$ and define $\delta_1,\ldots,\delta_n$ recursively by
\[
  \delta_1:=\varepsilon,\qquad
  \delta_{j+1}:=\Big(\frac{\delta_j}{3}\Big)^{2},
\]
so that $1\geq\delta_1\geq\delta_2\geq\cdots\geq\delta_n>0$ and
$3\sqrt{\delta_{j+1}}=\delta_j$ for every $j\in\{1,\ldots,n-1\}$. According
to \cref{lem:log} there exists some $x_\varepsilon\in\R$ such that
$\min_{k}p_k(x)\leq\delta_n$ holds for every $x\leq x_\varepsilon$. Fix
$x\leq x_\varepsilon$ and choose $k\in\{1,\ldots,n\}$ with
$p_k(x)\leq\delta_n\leq\delta_k$. If $k\geq 2$, then \cref{lem:cascade}
yields $p_{k-1}(x)\leq 3\sqrt{p_k(x)}\leq 3\sqrt{\delta_k}=\delta_{k-1}$.
Proceeding in the same manner we successively obtain $p_j(x)\leq\delta_j$
for $j=k,k-1,\ldots,1$ and, in particular,
$0\leq p_1(x)\leq\delta_1=\varepsilon$. Since $x\leq x_\varepsilon$ and
$\varepsilon\in(0,1]$ were arbitrary, this completes the proof.
\end{proof}

Having \cref{prop:top}, the strict case can now be settled by interchanging
the roles of $\gamma_1$ and $\gamma_2$.

\begin{theorem}\label{thm:strict}
Suppose that $s_1=s_2=0$. Then $M^*=m_*$ and $\gamma_1=\gamma_2^{(c)}$ holds
for $c:=M^*$.
\end{theorem}

\begin{proof}
All results established so far apply to the pair $(\gamma_2,\gamma_1)$ too,
since it fulfills the same assumptions as $(\gamma_1,\gamma_2)$. We will let
$\hat\sigma_k$, $\hat\varrho_k$, $\hat M$ and $\hat M^*$ denote the objects
corresponding to the pair $(\gamma_2,\gamma_1)$. By definition we have
$\hat\sigma_k=(\Pk{k}_1)^{-1}\circ\Pk{k}_2=\sigma_k^{-1}$ and hence
\begin{equation}\label{eq:swap}
  \hat\varrho_k\big(\sigma_k(s)\big)=\frac{s}{\sigma_k(s)}
  =\frac{1}{\varrho_k(s)}
\end{equation}
for all $s>0$ and $k\in\{1,\ldots,n\}$. Considering that each $\sigma_k$ is a
bijection of $(0,\infty)$, eq.~\eqref{eq:swap} implies
\[
  \hat M^*=\max_{k\in\{1,\ldots,n\}}\,\sup_{t>0}\hat\varrho_k(t)
  =\max_{k\in\{1,\ldots,n\}}\,\sup_{s>0}\frac{1}{\varrho_k(s)}
  =\frac{1}{m_*}.
\]
Applying \cref{prop:top} to both pairs yields
\[
  \lim_{s\downarrow 0}\varrho_1(s)=M^*\qquad\text{and}\qquad
  \lim_{t\downarrow 0}\hat\varrho_1(t)=\hat M^*=\frac{1}{m_*}.
\]
On the other hand, considering $\lim_{s\downarrow 0}\sigma_1(s)=0$ and using
eq.~\eqref{eq:swap} once more we get
\[
  \lim_{t\downarrow 0}\hat\varrho_1(t)
  =\lim_{s\downarrow 0}\hat\varrho_1\big(\sigma_1(s)\big)
  =\lim_{s\downarrow 0}\frac{1}{\varrho_1(s)}=\frac{1}{M^*}.
\]
Altogether we have shown $M^*=m_*$, i.e., $\sup_{s>0}M(s)=\inf_{s>0}m(s)$,
and applying \cref{cor:monotone} completes the proof.
\end{proof}

\subsection{Conclusion}\label{subsec:conclusion}

Combining \cref{prop:nonstrict} and \cref{thm:strict} yields the following
result, which is the announced reformulation of \cref{thm:main} in terms of
Williamson measures.

\begin{theorem}\label{thm:mainW}
Let $d\geq 2$ and $\gamma_1,\gamma_2\in\cP$. Then the Archimedean copulas
$C_{\gamma_1},C_{\gamma_2}\in\Car$ have the same Kendall distribution
function if, and only if, there exists some $c>0$ with
$\gamma_1=\gamma_2^{(c)}$.
\end{theorem}

\begin{proof}
If $\gamma_1=\gamma_2^{(c)}$ holds for some $c>0$, then
$C_{\gamma_1}=C_{\gamma_2}$ according to \cref{lem:scaling}. The reverse
implication has been proved in \cref{prop:nonstrict} (if $s_{\gamma_1}>0$ or
$s_{\gamma_2}>0$) and in \cref{thm:strict} (if
$s_{\gamma_1}=s_{\gamma_2}=0$).
\end{proof}

\begin{proof}[Proof of \cref{thm:main}]
According to \cref{thm:williamson} there exist $\gamma_1,\gamma_2\in\cP$ with
$C_i=C_{\gamma_i}$ for $i\in\{1,2\}$. \cref{thm:mainW} yields
$\gamma_1=\gamma_2^{(c)}$ for some $c>0$, and \cref{lem:scaling} implies
$C_1=C_{\gamma_2^{(c)}}=C_{\gamma_2}=C_2$.
\end{proof}

As a direct consequence of \cref{thm:mainW} we obtain injectivity of the
mapping assigning each suitably normalized Williamson measure the
corresponding Kendall distribution function:

\begin{corollary}\label{cor:injective}
Let $z_0\in(0,\infty)$ and $y_0\in(0,1)$. Then the mapping assigning each
$\gamma\in\cP$ with $\psi_\gamma(z_0)=y_0$ the Kendall distribution function
of $C_\gamma$ is injective.
\end{corollary}

\begin{proof}
Suppose that $\gamma_1,\gamma_2\in\cP$ fulfill
$\psi_{\gamma_1}(z_0)=\psi_{\gamma_2}(z_0)=y_0$ and induce the same Kendall
distribution function. According to \cref{thm:mainW} there exists some $c>0$
with $\gamma_1=\gamma_2^{(c)}$, hence
$\psi_{\gamma_2}(z_0/c)=\psi_{\gamma_1}(z_0)=y_0=\psi_{\gamma_2}(z_0)$ (see
the proof of \cref{lem:scaling}). Considering $y_0>0$ and using the fact that
$\psi_{\gamma_2}$ is strictly decreasing on
$\{z\geq 0\colon\psi_{\gamma_2}(z)>0\}$ we get $c=1$ and
$\gamma_1=\gamma_2$.
\end{proof}

\begin{remark}\label{rem:final}
We conclude with some comments on the proof of \cref{thm:mainW}:
\begin{enumerate}
\item The assumption that $\gamma_1$ and $\gamma_2$ induce the same Kendall
distribution function only enters via the identity $G_1=G_2\circ\sigma_n$ in
\cref{lem:cycle}, i.e., via the fact that $\varrho_1$ pursues $\varrho_n$,
which is essential for the monotonicity of $M$ and $m$. In contrast,
\cref{lem:sandwich}, the first assertion of \cref{lem:apriori},
eq.~\eqref{eq:pursuit} for $k\geq 2$, and \cref{lem:rate} hold for arbitrary
measures $\gamma_1,\gamma_2\in\cP$ with $s_{\gamma_1}=s_{\gamma_2}=0$.
\item For $d=2$ we have $n=1$ and $\sigma_0=\sigma_1$, so the numerator in
eq.~\eqref{eq:pursuit} vanishes identically. Accordingly, $M=m=\varrho_1$ is
non-increasing as well as non-decreasing, i.e., constant, and
\cref{cor:monotone} directly yields the bivariate result.
\item In the logarithmic coordinates used in \cref{subsec:cascade} the proof
can be summarized as follows: the functions $\tau_1,\ldots,\tau_n$ form a
pursuit system, which contracts as $x$ increases. In the non-strict case the
system starts in consensus at $x=\ln s_1$ and therefore stays in consensus.
In the strict case the system is defined on the whole real line and,
according to \cref{lem:apriori}, fulfills
$\max_k\tau_k-\min_k\tau_k\leq 2\ln n$. Since for $k\in\{2,\ldots,n\}$ a
function $\tau_k$ lying above its predecessor $\tau_{k-1}$ is pulled down at
a rate which is bounded from below (\cref{lem:rate}), the supremum of the
system can only be approached for $x\to-\infty$ if $\tau_1$ approaches it
too (\cref{prop:top}); by symmetry the same holds for the infimum, which is
only possible if the system is in consensus.
\end{enumerate}
\end{remark}

\section*{Declaration on the use of AI tools}
During the development of this work we used the AI tool Claude (Anthropic) to
explore proof strategies and to revise the
language of the manuscript. We have checked all statements and proofs and take
full responsibility for the content of this paper.




\begin{thebibliography}{9}

\bibitem{DS}
F.~Durante and C.~Sempi.
\newblock {\em Principles of copula theory}.
\newblock Taylor \& Francis Group LLC, Boca Raton FL, 2016.

\bibitem{GNZ}
C.~Genest, J.~Ne\v{s}lehov\'a, and J.~Ziegel.
\newblock Inference in multivariate {A}rchimedean copula models.
\newblock {\em TEST}, 20(2):223--256, 2011.

\bibitem{GR}
C.~Genest and L.-P.~Rivest.
\newblock Statistical inference procedures for bivariate {A}rchimedean
  copulas.
\newblock {\em Journal of the American Statistical Association},
  88(423):1034--1043, 1993.

\bibitem{KDT}
T.~M.~Kasper, N.~Dietrich, and W.~Trutschnig.
\newblock On convergence and mass distributions of multivariate {A}rchimedean
  copulas and their interplay with the {W}illiamson transform.
\newblock {\em Journal of Mathematical Analysis and Applications},
  529(1):127555, 2024.
\newblock (Theorem and equation numbers refer to the preprint
  arXiv:2210.11868v1.)

\bibitem{MN}
A.~J.~McNeil and J.~Ne\v{s}lehov\'a.
\newblock Multivariate {A}rchimedean copulas, $d$-monotone functions and
  $\ell_1$-norm symmetric distributions.
\newblock {\em The Annals of Statistics}, 37(5B):3059--3097, 2009.

\bibitem{Nelsen}
R.~B.~Nelsen.
\newblock {\em An Introduction to Copulas}.
\newblock Springer-Verlag, Berlin Heidelberg, 2006.

\bibitem{Williamson}
R.~E.~Williamson.
\newblock Multiply monotone functions and their {L}aplace transforms.
\newblock {\em Duke Mathematical Journal}, 23(2):189--207, 1956.

\end{thebibliography}
\end{document}